\documentclass[a4paper]{amsart}

\usepackage{ucs}
\usepackage[utf8x]{inputenc}
\usepackage[ngerman, english]{babel}

\usepackage{babelbib}

\usepackage{mathtools}

\usepackage{bbm}
\usepackage{mathrsfs}

\usepackage{amsmath, amssymb, amsthm}
\usepackage{IEEEtrantools}

\makeatletter
\def\blfootnote{\xdef\@thefnmark{}\@footnotetext}
\makeatother

\newcommand{\abs}[1]{\left| #1 \right|}
\newcommand{\BLOBS}{\mathscr{L}}
\newcommand{\Cb}{C_{\text{b}}}
\newcommand{\CN}{\mathbb{C}}
\newcommand{\D}{\mathrm{d}}
\newcommand{\dup}[2]{\left\langle #1, #2 \right\rangle}
\newcommand{\FT}{\mathcal{F}}
\newcommand{\Ind}{\mathbbmss{1}}
\newcommand{\iu}{\mathrm{i}}
\newcommand{\iv}[1]{\frac{1}{ #1}}
\newcommand{\jb}[1]{\langle #1 \rangle}
\newcommand{\LSS}{\mathcal{S}}
\newcommand{\N}{{\mathbb{N}}}
\newcommand{\norm}[1]{\left\Vert #1 \right\Vert}
\newcommand{\R}{{\mathbb{R}}}
\newcommand{\set}[1]{\left\{ #1 \right\}}
\DeclareMathOperator{\supp}{supp}
\newcommand{\Z}{{\mathbb Z}}

\theoremstyle{plain}
\newtheorem{theorem}{Theorem}
\newtheorem{lemma}[theorem]{Lemma}

\theoremstyle{definition}

\usepackage{hyperref}

\title[LWP for the NLS in $M_{p, q}^s \cap M_{\infty, 1}$]
{Local well-posedness for the nonlinear Schr\"odinger equation
in the intersection of modulation spaces
$M_{p, q}^s(\R^d) \cap M_{\infty, 1}(\R^d)$}

\author{L. Chaichenets}
\address{Leonid Chaichenets, Department of Mathematics,
Institute for Analysis, Karlsruhe Institute of Technology,
76128 Karlsruhe, Germany}
\email{leonid.chaichenets@kit.edu}
\author{D. Hundertmark}
\address{Dirk Hundertmark, Department of Mathematics,
Institute for Analysis, Karlsruhe Institute of Technology,
76128 Karlsruhe, Germany}
\email{dirk.hundertmark@kit.edu}
\author{P. Kunstmann}
\address{Peer Christian Kunstmann, Department of Mathematics,
Institute for Analysis, Karlsruhe Institute of Technology,
76128 Karlsruhe, Germany}
\email{peer.kunstmann@kit.edu}
\author{N. Pattakos}
\address{Nikolaos Pattakos, Department of Mathematics,
Institute for Analysis, Karlsruhe Institute of Technology,
76128 Karlsruhe, Germany}
\email{nikolaos.pattakos@kit.edu}

\begin{document}
\blfootnote{\copyright 2019 by the authors. Faithful reproduction of
this article, in its entirety, by any means is permitted for
noncommercial purposes.}
\maketitle

\begin{abstract}
We introduce a Littlewood-Paley characterization of modulation
spaces and use it to give an alternative proof of the algebra property,
somehow implicitly contained in \cite{sugimoto2011}, of the intersection
$M^s_{p,q}(\R^d) \cap M_{\infty, 1}(\R^d)$ for
$d \in \N$, $p, q \in [1, \infty]$ and $s \geq 0$.
We employ this algebra property to show the local well-posedness of the
Cauchy problem for the cubic nonlinear Schr\"odinger equation in the above
intersection. This improves \cite[Theorem 1.1]{benyi2009} by B\'enyi and
Okoudjou, where only the case $q = 1$ is considered, and closes a gap in
the literature. If $q > 1$ and $s > d \left(1 - \frac{1}{q}\right)$ or
if $q = 1$ and $s \geq 0$ then
$M^s_{p,q}(\R^d) \hookrightarrow M_{\infty, 1}(\R^d)$ and the above
intersection is superfluous. For this case we also reobtain a
H\"older-type inequality for modulation spaces.
\end{abstract}

\section{Introduction}
In this paper we contribute to the general theory of modulation spaces.
Modulation spaces $M_{p,q}^s(\R^d)$ were introduced by Feichtinger in
\cite{feichtinger1983}. Here, we only briefly recall their definition
and refer to Section \ref{sec:preliminaries} and the literature
mentioned there for more information. Fix a so-called
\emph{window function} $g \in \LSS(\R^d) \setminus \set{0}$. The
\emph{short-time Fourier transform} $V_gf$ of a tempered distribution
$f \in \LSS'(\R^d)$ with respect to the window $g$ is defined by
\begin{equation}
\label{eqn:STFT}
(V_gf)(x, \xi) = \frac{1}{(2 \pi)^{\frac{d}{2}}}
\overline{\dup{f}{M_{\xi} S_x g}} \qquad \forall x, \xi \in \R^d,
\end{equation}
where $S_x g(y) = g(y - x)$ denotes the \emph{right-shift} by
$x \in \R^d$, $\left(M_\xi g\right)(y) = e^{\iu k \cdot y} g(y)$
the \emph{modulation} by $\xi \in \R^d$ and $\dup{f}{g} =
\int_{\R^d} \overline{f}(x) g(x) \D{x}$ for
$f \in L^1_{\text{loc}}(\R^d)$, $g \in \LSS(\R)^d$. We define
\begin{eqnarray*}
M_{p,q}^s(\R^d) & = &
\set{f \in \LSS'(\R^d) \Big| \norm{f}_{M_{p,q}^s(\R^d)} < \infty}
\text{, where} \\
\norm{f}_{M_{p,q}^s(\R^d)} & = & \norm{\xi \mapsto
\jb{\xi}^s\norm{V_gf\left(\cdot, \xi\right)}_p}_q
\end{eqnarray*}
for $s \in \R$, $p, q \in [1, \infty]$. As usual in the literature, we
set $M_{p, q}(\R^d) \coloneqq M_{p,q}^0(\R^d)$ and often shorten the
notation for $M_{p,q}^s(\R^d)$ to $M_{p,q}^s$. It can be shown, that
the $M_{p,q}^s(\R^d)$ are Banach spaces and that different choices of
the window function $g$ lead to equivalent norms.

To state our first result, let us recall the definition of the
Littlewood-Paley decomposition. Consider a smooth radial function
$\phi_0 \in C_c^\infty(\R^d)$ with $\phi_0(\xi) = 1$ for all
$\abs{\xi} \leq \frac{1}{2}$ and $\supp(\phi_0) \subseteq B_1(0)$. Set
$\phi_1 = \phi_0\left(\frac{\cdot}{2}\right) - \phi_0$ and
$\phi_l = \phi_1\left(\frac{\cdot}{2^{l-1}}\right)$ for all $l \in \N$.
The multiplier operators defined by
\begin{equation*}
\Delta_l f \coloneqq
\frac{1}{(2\pi)^{\frac{d}{2}}} \check{\phi}_l \ast f =
\FT^{(-1)} \phi_l \FT f \qquad
\forall l \in \N_0 \, \forall f \in \LSS'(\R^d)
\end{equation*}
are called \emph{dyadic decomposition operators} and the sequence
$(\Delta_l f)_{l \in \N_0}$ is called the \emph{Littlewood-Paley
decomposition} of $f \in \LSS'(\R^d)$. Above, $\FT$ denotes the
usual \emph{Fourier transform} and $\FT^{(-1)}$ its inverse.

Our first result is
\begin{theorem}[Littlewood-Paley characterization]
\label{thm:littlewood_paley}
Let $d \in \N$, $p,q \in [1, \infty]$ and $s \in \R$. Then
\begin{equation*}
\norm{f} \coloneqq
\norm{\left(2^{ls}
\norm{\Delta_l f}_{M_{p,q}(\R^d)}\right)_{l \in \N_0}}_q
\qquad \forall f \in \LSS'(\R^d)
\end{equation*}
defines an equivalent norm on $M_{p,q}^s(\R^d)$. The constants of the norm
equivalence depend only on $d$ and $s$.
\end{theorem}

The above characterization of modulation spaces is new and we shall use
it to prove that the intersections
$M_{p,q}^s(\R^d) \cap M_{\infty, 1}(\R^d)$ are
\emph{Banach *-algebras\footnotemark}.
\footnotetext{For us, a Banach *-algebra $X$ is a Banach algebra over
$\CN$ on which a continuous \emph{involution} $*$ is defined, i.e.
$(x+y)^* = x^* + y^*$, $(\lambda x)^* = \overline{\lambda} x^*$,
$(xy)^* = y^* x^*$ and $(x^*)^* = x$ for any $x,y \in X$ and
$\lambda \in \CN$. We neither require $X$ to have a unit nor $C = 1$ in
the estimates $\norm{x \cdot y} \leq C \norm{x} \norm{y}$,
$\norm{x^*} \leq C \norm{x}$.} To state this second result, let us
denote by $\Cb(\R^d)$ the space of bounded complex-valued continuous
functions on $\R^d$, where $d \in \N$. We then have
\begin{theorem}[Algebra property]
\label{thm:algebra}
Let $d \in \N$, $p, q \in [1, \infty]$ and $s \geq 0$. Then
$M_{p,q}^s(\R^d) \cap M_{\infty, 1}(\R^d)$ is a Banach *-algebra with
respect to pointwise multiplication and complex conjugation. These
operations are well-defined due to the embedding
$M_{\infty, 1}(\R^d) \hookrightarrow \Cb(\R^d)$ Furthermore, if $q > 1$
and $s > d \left(1 - \frac{1}{q}\right)$ or if $q = 1$, then
$M^s_{p,q}(\R^d) \hookrightarrow M_{\infty, 1}(\R^d)$, so in particular
$M_{p,q}^s(\R^d)$ is a Banach *-algebra, in that case.
\end{theorem}
The latter case of Theorem \ref{thm:algebra} had been
observed already in 1983 by Feichtinger in his aforementioned pioneering
work on modulation spaces (cf. \cite[Proposition 6.9]{feichtinger1983}),
where he proves it using a rather abstract approach via Banach
convolution triples. The case $q > 1$ and
$s \in \left[0, d \left(1 - \frac{1}{q}\right)\right]$ seems to be new,
at least as a statement. A different proof of Theorem \ref{thm:algebra}
can be given following the idea of proof of
\cite[Proposition 3.2]{sugimoto2011}, see
\cite[Proposition 4.2]{chaichenets2018}.

Our third result is a H\"older-type inequality for modulation spaces,
which is stated in
\begin{theorem}[H\"older-type inequality]
\label{thm:hoelder}
Let $d \in \N$ and $p, p_1, p_2, q \in [1, \infty]$ be such that
$\iv{p} = \iv{p_1} + \iv{p_2}$. For $q > 1$ let
$s > d\left(1 - \frac{1}{q}\right)$ and for $q = 1$ let $s \geq 0$.
Then there is a $C > 0$ such that for any $f \in M_{p_1,q}^s(\R^d)$ and
any $g \in M_{p_2,q}^s(\R^d)$ one has $f g \in M_{p, q}^s(\R^d)$ and
\begin{equation}
\label{eqn:hoelder}
\norm{fg}_{M_{p,q}^s(\R^d)} \leq C
\norm{f}_{M_{p_1,q}^s(\R^d)} \norm{g}_{M_{p_2,q}^s(\R^d)}.
\end{equation} 
\end{theorem}
The above pointwise multiplication $fg$ is well-defined due to the
embedding formulated in Theorem \ref{thm:algebra}. The
constant $C$ does \emph{not} depend on $p$, $p_1$ or $p_2$.

Theorem \ref{thm:hoelder} easily generalizes to $m \in \N$ factors and
$p, p_1, \ldots, p_m \in (0, \infty]$. Hence, it extends the
multilinear estimate from \cite[Equation 2.4]{benyi2009} to the case
$q_0 = \ldots = q_m > 1$. The case $q \neq \infty$ and $p \neq \infty$
of our Theorem \ref{thm:hoelder} is a special case of
\cite[Theorem 1.4]{guo2018}. That theorem also establishes that the
condition $q = 1$ and $s \geq 0$ or $q > 1$ and $s > \frac{d}{q'}$ is
necessary for the space $M_{p, q}^s(\R^d)$ to be a Banach *-algebra.

Let us also mention, that other results on H\"older-type inequalities in
modulation spaces, i.e.,
\begin{equation*}
\norm{f g}_{M_{p, q}^s} \lesssim
\norm{f}_{M_{p_1, q_1}^{s_1}} \norm{g}_{M_{p_2, q_2}^{s_2}},
\end{equation*}
include \cite[Theorem 2.4]{toft2015} and
\cite[Proposition 3.5]{cordero2009} (but see also
\cite[Lemma 4.1]{wang2006}). In the case $q = q_1 = q_2$, the
former result requires $q \leq 2$ and $p \neq \infty$, while the latter
result only covers the case $q = 1$.

Here we present a direct proof of Theorem \ref{thm:hoelder}, close to
the approach found in \cite[Corollary 4.2]{wang2006} and involving an
application of Theorem \ref{thm:algebra}. For a proof avoiding
the Littewood-Paley characterization see the proof of
\cite[Theorem 4.3]{chaichenets2018}. Yet another and more abstract
proof could be given by invoking \cite[Theorem 3]{feichtinger1980} for a
specific choice of Banach convolution triples.

Lastly, we employ Theorem \ref{thm:algebra} to study the Cauchy problem
for the cubic nonlinear Schr\"odinger equation
(\emph{NLS})
\begin{equation}
\label{eqn:Cauchy_NLS}
\left\{
\begin{IEEEeqnarraybox}[][c]{rCl?rCl}
\iu \frac{\partial u}{\partial t} (x, t) +
\Delta u(x,t) \pm
\abs{u}^2 u(x, t) & = & 0 &
(x,t) & \in & \R^d \times \R, \\
u(x, 0) & = & u_0(x) & x & \in & \R^d,
\end{IEEEeqnarraybox}
\right.
\end{equation}
where the initial data $u_0$ is in an intersection of modulation spaces
$M_{p,q}^s(\R^d) \cap M_{\infty, 1}(\R^d)$. We are interested in
\emph{mild solutions} $u$ of \eqref{eqn:Cauchy_NLS}, i.e.
\begin{equation*}
u \in C\left([0,T), M_{p,q}^s(\R^d) \cap M_{\infty, 1}(\R^d) \right)
\end{equation*}
for some $T > 0$ which satisfy the corresponding integral equation
\begin{equation}
\label{eqn:duhamel}
u(\cdot, t) =
e^{\iu t \Delta} u_0 \pm
\iu \int_0^{t} e^{\iu(t - \tau) \Delta}
\left(\abs{u}^2u(\cdot, \tau)\right) \D{\tau} \qquad
\forall t \in [0, T).
\end{equation}
Our last result is stated in 
\begin{theorem}[Local well-posedness]
\label{thm:lwp}
Let $d \in \N$, $p \in [1,\infty]$, $q \in [1, \infty)$ and
$s \geq 0$. Set
$X = M_{p,q}^s(\R^d) \cap M_{\infty, 1}(\R^d)$ and $X(T) = C([0,T], X)$,
$X_*(T) = C([0,T), X)$ for any $T > 0$. Assume
that $u_0 \in X$. Then, there exists a unique maximal mild solution
$u \in X_*(T_*)$ of \eqref{eqn:Cauchy_NLS} and the
\emph{blow-up alternative}
\begin{equation*}
T_* < \infty \qquad \Rightarrow \qquad
\limsup_{t \to T_*-} \norm{u(\cdot, t)}_X = \infty
\end{equation*}
holds. Moreover, for any $T' \in (0, T_*)$ there exists a neighborhood
$V$ of $u_0$ in $X$, such that the initial-data-to-solution-map
$V \to X(T')$, $v_0 \mapsto v$ is Lipschitz continuous.
\end{theorem}

As already stated in Theorem \ref{thm:algebra} one has that, if
$q > 1$ and $s > d \left(1 - \frac{1}{q}\right)$ or if
$q = 1$, then $M^s_{p,q}(\R^d) \hookrightarrow M_{\infty, 1}(\R^d)$ and
so $X = M_{p,q}^s(\R^d)$, in that case.

In the case $q = \infty$ excluded in Theorem \ref{thm:lwp}, the
situation is more subtle. Following our proof, one obtains local
well-posedness in the larger space
\begin{equation*}
L^\infty([0, T), M_{p, \infty}^s(\R^d) \cap M_{\infty, 1}(\R^d)).
\end{equation*}
The missing continuity in time is due to the properties of the free
Schr\"odinger evolution and we refer to the remarks after Theorem
\ref{thm:schroedinger}.

The precursors of Theorem \ref{thm:lwp} are \cite[Theorem 1.1]{wang2006}
by Wang, Zhao and Guo for the space $M_{2,1}^0(\R^d)$ and
\cite[Theorem 1.1]{benyi2009} due to B\'enyi and Okoudjou for the space
$M_{p,1}^s(\R^d)$ with $p \in [1, \infty]$ and $s \geq 0$. In fact,
Theorem \ref{thm:lwp} generalizes \cite[Theorem 1.1]{benyi2009} to
$q \geq 1$: Although our theorem is stated for the cubic nonlinearity,
this is for simplicity of the presentation only. The proof allows for an
easy generalization to \emph{algebraic nonlinearities} considered in
\cite{benyi2009}, which are of the form
\begin{equation}
\label{eqn:algebraic_nonlinearity}
f(u) = g(\abs{u}^2)u = \sum_{k=0}^\infty c_k \abs{u}^{2k}u,
\end{equation}
where $g$ is an entire function. Also,
\cite[Theorems 1.2 and 1.3]{benyi2009}, which concern the nonlinear wave
and the nonlinear Klein-Gordon equation respectively, can be generalized
in the same spirit. The reason for this is that the proof of these
results is based on the well-known Banach's contraction principle, on
the fact that the free propagator is a $C_0$-group, and on the algebra
property of the spaces under consideration. Although the ingredients
seem to be known in the community, the results to be found in the
literature (e.g. \cite[Theorem 6.2]{wang2011}) are not as general as
Theorem \ref{thm:lwp}. In this sense, it closes a gap in the literature.

Let us remark that local well-posedness results in the case of
modulation spaces that are not Banach *-algebras are
\cite[Theorem 1.4]{guo2016} for $u_0 \in M_{2, q}(\R)$ with
$q \in [2, \infty)$ and
\cite[Theorem 6]{pattakos2018a} with $u_0 \in M_{p, q}^s(\R)$ with
either
$p \in [2, 3]$, $q \in \left[1, \frac{3}{2} \right]$ and $s \geq 0$ or
$p \in [2, 3]$, $q \in \left(\frac{3}{2}, \frac{18}{11}\right]$ and
$s > \frac{2}{3} - \iv{q}$ or
$q \in \left(\frac{18}{11}, 2\right]$,
$p \in \left[2, \frac{10 q}{7q - 6}\right)$ and $s > \frac{2}{3} - \iv{q}$
(see also \cite[Theorem 4]{pattakos2018}).

The remainder of our paper is structured as follows. We start with
Section \ref{sec:preliminaries} providing an overview over modulation
spaces and the free Schr\"odinger propagator on them. In Section
\ref{sec:littlewood_paley} we apply methods from the Littlewood-Paley
theory to prove Theorem \ref{thm:littlewood_paley}. In
the subsequent Section \ref{sec:algebra} we prove the algebra property
from Theorem \ref{thm:algebra}, notice the resulting similar
property for weighted sequence spaces in Lemma
\ref{lem:algebra_property_seq}, and deduce the H\"older-type inequality
stated in Theorem \ref{thm:hoelder}. Finally, we prove Theorem
\ref{thm:lwp} on the local well-posedness in Section
\ref{sec:local_well-posedness}.

\subsection*{Notation}
We denote generic constants by $C$. To emphasize on which quantities
a constant depends we write e.g. $C = C(d)$ or $C = C(d,s)$. Sometimes
we omit a positive constant from an inequality by writing “$\lesssim$”,
e.g. $A \lesssim_d B$ instead of $A \leq C(d) B$. By $A \approx B$ we
mean $A \lesssim B$ and $B \lesssim A$. Special constants are $d \in \N$
for the \emph{dimension}, $p, q \in [1,\infty]$ for the \emph{Lebesgue}
exponents and $s \in \R$ for the \emph{regularity} exponent. By $p'$ we
mean the \emph{dual} exponent of $p$, that is the number satisfying
$\frac{1}{p} + \frac{1}{p'} = 1$.

We denote by $\LSS(\R^d)$ the set of \emph{Schwartz functions} and by
$\LSS'(\R^d)$ the space of \emph{tempered distributions}. Furthermore,
we denote the \emph{Bessel potential spaces} or simply $L^2$-based
\emph{Sobolev spaces} by $H^s = H^s(\R^d)$. For the space of smooth
functions with compact support we write $C_c^\infty$. The letters
$f, g, h$ denote either generic functions $\R^d \to \CN$ or generic
tempered distributions and $(a_{k})_{k \in \Z^d} =
(a_{k})_{k} = (a_k)$, $(b_{k})_{k \in \Z^d} = (b_{k})_{k} = (b_{k})$
denote generic complex-valued sequences. By $\jb{\cdot} =
\sqrt{1 + \abs{\cdot}^2}$ we mean the \emph{Japanese bracket}.

For a Banach space $X$ we write $X^*$ for its dual and $\norm{\cdot}_X$
for the norm it is canonically equipped with. By $\BLOBS(X,Y)$ we denote
the space of all bounded linear maps from $X$ to $Y$, where $Y$ is
another Banach space, and set $\BLOBS(X) = \BLOBS(X,X)$. By
$[X,Y]_\theta$ we mean complex interpolation between $X$ and $Y$, if
$(X,Y)$ is an interpolation couple. For brevity we write
$\norm{\cdot}_p$ for the $p$-norm on the \emph{Lebesgue space}
$L^p = L^p(\R^d)$, the \emph{sequence space}
$l^p = l^p(\Z^d)$ or $l^p = l^p(\N_0)$ and $\norm{(a_k)}_{q,s} \coloneqq
\norm{\left(\jb{k}^sa_k\right)}_q$ for the norm on
$\jb{\cdot}^s$-weighted sequence spaces $l_s^q = l_s^q(\Z^d)$. If the
norm is apparent from the context, we write $B_r(x)$ for a ball of
radius $r$ around $x \in X$.

We use the symmetric choice of constants for the Fourier transform and
also write
\begin{eqnarray*}
\hat{f}(\xi) & \coloneqq &
(\FT f)(\xi) =
\frac{1}{(2 \pi)^{\frac{d}{2}}}
\int_{\R^d} e^{-\iu \xi \cdot x} f(x) \D{x}, \\
\check{g}(x) & \coloneqq & \left(\FT^{(-1)} g\right)(x) =
\frac{1}{(2 \pi)^{\frac{d}{2}}}
\int_{\R^d} e^{\iu \xi \cdot x} g(\xi) \D{\xi}.
\end{eqnarray*}

\section{Preliminaries}
\label{sec:preliminaries}
As already mentioned in the introduction, modulation spaces were
introduced by Feichtinger in \cite{feichtinger1983} in the setting of
locally compact Abelian groups. A thorough introduction is given in the
textbook \cite{groechenig2001} by Gr\"ochenig. A presentation
incorporating the characterization of modulation spaces via
\emph{isometric decomposition operators}, which we are going to use
below, is contained in the paper \cite[Section 2, 3]{wang2007} by Wang
and Hudzik. A survey on modulation spaces and nonlinear evolution
equations is given in \cite{ruzhansky2012}.

A convenient equivalent norm on modulation spaces which we are
going to use is constructed as follows
(cf. \cite[Propostition 2.1]{wang2007}): Set
$Q_{0} \coloneqq \left[-\frac{1}{2}, \frac{1}{2} \right)^d$ and
$Q_{k} \coloneqq Q_{0} + k$ for all $k \in \Z^d$.
Consider a smooth partition of unity
$(\sigma_{k})_{k \in \Z^d} \in
\left(C_c^\infty(\R^d)\right)^{\Z^d}$ satisfying
\begin{itemize}
\item
$\exists c > 0: \, \forall k \in \Z^d: \,
\forall \eta \in Q_{k}: \,
\abs{\sigma_{k}(\eta)} \geq c$,
\item
$\forall k \in \Z^d: \,
\supp(\sigma_{k}) \subseteq B_{\sqrt{d}}\left(k\right)$,
\item
$\sum_{k \in \Z^d} \sigma_{k} = 1$,
\item
$\forall m \in \N_0: \, \exists C_m > 0: \,
\forall k \in \Z^d: \,
\forall \alpha \in \N_0^d: \,
\abs{\alpha} \leq m \Rightarrow
\norm{D^\alpha \sigma_{k}}_\infty \leq C_m$
\end{itemize}
and define the \emph{isometric decomposition operators}
$\Box_{k} \coloneqq \FT^{(-1)} \sigma_{k} \FT$. We need the following often
used (cf. \cite[Proposition 1.9]{wang2007})
\begin{lemma}[Bernstein multiplier estimate]
\label{lem:bernstein}
Let $d \in \N$, $\sigma \in \LSS(\R^d)$ and
$r, p_1, p_2 \in [1, \infty]$ such that
$1 + \frac{1}{p_2} = \frac{1}{r} + \frac{1}{p_1}$.
Consider the \emph{multiplier operator}
$T_\sigma: \LSS'(\R^d) \to \LSS'(\R^d)$ with \emph{symbol} $\sigma$
defined by
\begin{equation*}
T_\sigma f = \FT^{(-1)} \sigma \FT f =
\frac{1}{(2\pi)^{\frac{d}{2}}} \check{\sigma} \ast f \qquad
\forall f \in \LSS'(\R^d).
\end{equation*}
Then, for any $f \in \LSS'(\R^d)$, every derivative of
$T_\sigma f \in C^\infty(\R^d)$ (including $T_\sigma f$) has at most
polynomial growth. Furthermore $\norm{T_\sigma f}_{p_2} \leq
\frac{\norm{\hat{\sigma}}_r}{(2\pi)^{\frac{d}{2}}} \norm{f}_{p_1}$ for
any $f \in L^{p_1}(\R^d)$.
\end{lemma}
Putting $r = 1$ and $p_1 = p_2 = p$ in Lemma \ref{lem:bernstein}, shows
that $\Box_{k} f \in C^\infty(\R^d)$ for $f \in \LSS'(\R^d)$ and
$\norm{\Box_{k}}_{\BLOBS(L^p(\R^d))}$ is bounded independently of $k$
and $p$. The aforementioned equivalent norm for the modulation space
$M_{p,q}^s(\R^d)$ is given by (see \cite[Proposition 2.1]{wang2007})
\begin{equation}
\label{eqn:equivalent_norm}
\norm{f}_{M_{p,q}^s} \approx
\norm{\left( \jb{k}^s
\norm{\Box_{k}f}_p
\right)_{k \in \Z^d}}_q.
\end{equation}
Choosing a different partition of unity $(\sigma_{k})$ yields yet
another equivalent norm.

\begin{lemma}[Continuous embeddings]
Let $s_1 \geq s_2$, $1 \leq p_1 \leq p_2 \leq \infty$,
$1 \leq q_1 \leq q_2 \leq \infty$, $q > 1$ and $s > \frac{d}{q'}$. Then
\label{lem:embeddings}
\begin{enumerate}
\item
\label{it:interembeddings}
$M_{p_1, q_1}^{s_1}(\R^d) \subseteq M_{p_2, q_2}^{s_2}(\R^d)$
and the embedding is continuous,

\item
\label{it:regularity_fourier_integrability}
$M_{p_1,q}^s(\R^d) \subseteq M_{p_1,1}(\R^d)$ and the embedding is
continuous,

\item
\label{it:contembedding}
$M_{p_1,1}(\R^d) \hookrightarrow C_b(\R^d)$.
\end{enumerate}
\end{lemma}
Lemma \ref{lem:embeddings} is well-known (cf.
\cite[Proposition 2.5, 2.7]{wang2007}), but for convenience we sketch a
\begin{proof}
\begin{enumerate}
\item
One can change indices one by one. The inclusion for ``$s$'' is by
monotonicity and the inclusion for ``$q$'' is by the embeddings of the
$l^q$ spaces. For the ``$p$''-embedding consider
$\tau \in C_c^\infty(\R^d)$ such that $\tau|_{B_{\sqrt{d}}} \equiv 1$
and $\supp(\tau) \subseteq B_d$. For every $k \in \Z^d$, consider the
shifted symbol $\tau_{k} = S_k \tau$, define the corresponding
multiplier operator $\tilde{\Box}_{k} = \FT^{(-1)} \tau_{k} \FT$ and
observe, that $\hat{\tau}_k = M_k \hat{\tau}$. Hence, by Lemma
\ref{lem:bernstein}, the family
$\left(\tilde{\Box}_{k}\right)_{k \in \Z^d}$ is bounded in
$\BLOBS(L^{p_1}(\R^d),L^{p_2}(\R^d))$. So,
$\norm{\Box_{k} f}_{p_2} =
\norm{\tilde{\Box}_{k} \Box_{k} f}_{p_2} \lesssim_d
\norm{\Box_{k}f}_{p_1}$ for any $k \in \Z^d$. Recalling
\eqref{eqn:equivalent_norm} completes the argument.

\item By H\"older's inequality we immediately have
\begin{eqnarray*}
\norm{f}_{p_1,1} & \approx &
\sum_{k \in \Z^d} \norm{\Box_k f}_{p_1} \leq
\left(\sum_{k \in \Z^d} \jb{k}^{-sq'}\right)^{\frac{1}{q'}}
\left(\sum_{k \in \Z^d} \jb{k}^{sq}
\norm{\Box_k f}_p^q \right)^{\frac{1}{q}} \\
& \approx &
\left(\sum_{k \in \Z^d} \jb{k}^{-sq'}\right)^{\frac{1}{q'}}
\norm{f}_{M_{p_1,q}^s}
\end{eqnarray*}
and the first factor is finite for $s > \frac{d}{q'}$ by comparison with
the integral $\int_{\R^d} \jb{x}^{-sq'} \D{x}$.

\item By part (\ref{it:interembeddings}) it is enough to show that
$M_{\infty,1} \hookrightarrow C_b$. For any $f \in M_{\infty,1}$
we have $\underbrace{\sum_{\abs{k} \leq N} \Box_{k} f}_{\in C^\infty}
\to f$ in $\LSS'$ as $N \to \infty$. But simultaneously,
the series $\sum_{k \in \Z^d} \Box_{k} f$ is absolutely convergent in
$L^\infty$ to, say, $g \in C_b$. As $M_{\infty, 1} \hookrightarrow S'$
(see \cite[Thm. 6.1 (B)]{feichtinger1983}), we have $f = g$.
\end{enumerate}
\end{proof}

For the proof of Theorem \ref{thm:algebra} we will need
the following (cf. \cite[eqn. (2.4)]{benyi2009}).
\begin{lemma}[Bilinear estimate]
\label{lem:bilinear}
Let $d \in \N$ and $1 \leq p \leq \infty$. Assume $f \in M_{p, q}(\R^d)$
and $g \in M_{\infty, 1}(\R^d)$. Then
\begin{equation*}
\norm{fg}_{M_{p,q}(\R^d)} \lesssim
\norm{f}_{M_{p, q}(\R^d)}
\norm{g}_{M_{\infty, 1}(\R^d)},
\end{equation*}
where the implicit constant does \emph{not} depend on $p$ or $q$.
\end{lemma}
For convenience, and because we will generalize Lemma \ref{lem:bilinear}
to Theorem \ref{thm:hoelder}, we present a proof close to the one of
\cite[Corollary 4.2]{wang2006}.
\begin{proof}
We use \eqref{eqn:equivalent_norm} to estimate the modulation space
norm of the left-hand side. Fix a $k \in \Z^d$. By the definition of
the operator $\Box_{k}$ we have
\begin{equation*}
\Box_{k}(fg) =
\frac{1}{(2\pi)^{\frac{d}{2}}} \FT^{(-1)}
\left(\sigma_{k} (\hat{f} \ast \hat{g}) \right) =
\frac{1}{(2\pi)^{\frac{d}{2}}}
\sum_{l,m \in \Z^d} \FT^{(-1)}
\left(\sigma_{k} ((\sigma_{l} \hat{f}) \ast
(\sigma_{m} \hat{g})) \right).
\end{equation*}
As the supports of the partition of unity are compact, many
summands vanish. Indeed, for any
$k, l, m \in \Z^d$
\begin{eqnarray*}
\supp\left(\sigma_{k} \left((\sigma_{l} \hat{f}) \ast
(\sigma_{m} \hat{g}) \right) \right) & \subseteq &
\supp(\sigma_{k}) \cap
\left(\supp(\sigma_{l}) + \supp(\sigma_{m})\right) \\
& \subseteq &
B_{\sqrt{d}}(k) \cap B_{2 \sqrt{d}}(l + m)
\end{eqnarray*}
and so $\sigma_{k} \left((\sigma_{l} \hat{f}) \ast
(\sigma_{m} \hat{g}) \right) \equiv 0$ if
$\abs{(k-l) - m} > 3\sqrt{d}$. Hence,
the double series over $l, m \in \Z^d$ boils down to a finite
sum of discrete convolutions
\begin{eqnarray*}
\Box_{k}(fg) & = &
\frac{1}{(2\pi)^{\frac{d}{2}}} \FT^{(-1)}
\left(\sigma_{k} \sum_{m \in M} \sum_{l \in \Z^d}
(\sigma_{l} \hat{f}) \ast
(\sigma_{k - l + m} \hat{g}) \right) \\
& = &
\Box_{k}
\sum_{m \in M} \sum_{l \in \Z^d}
(\Box_{l} f) \cdot (\Box_{k + m - l} g),
\end{eqnarray*}
where $M = \set{m \in \Z^d \Big| \, \abs{m} \leq 3 \sqrt{d}}$ and
$\# M \leq \left(6 \sqrt{d} + 1 \right)^d < \infty$. That was the
job of $\Box_{k}$ and we now get rid of it,
\begin{equation*}
\norm{\Box_{k}(fg)}_p \lesssim 
\sum_{m \in M} \sum_{l \in \Z^d}
\norm{(\Box_{l} f) \cdot 
(\Box_{k + m - l} g)}_p,
\end{equation*}
using the Bernstein multiplier estimate from Lemma \ref{lem:bernstein}.

Invoking H\"older's inequality we further estimate
\begin{equation}
\label{eqn:gen_hoelder}
\norm{\Box_{k}(fg)}_p \lesssim
\sum_{m \in M} \left(
\left(\norm{\Box_{l}(f)}_p\right)_{l} \ast
\left(\norm{\Box_{n + m}(g)}_{\infty}\right)_{n}\right)(k)
\end{equation}
pointwise in $k \in \Z^d$, where $\ast$ denotes the convolution of
sequences, and hence obtain
\begin{equation*}
\norm{fg}_{M_{p,q}} \lesssim
\norm{\left(\norm{\Box_{l} f}_p \right)_{l}}_q
\norm{\left(\norm{\Box_{n} g}_{\infty} \right)_n}_1
\end{equation*}
by Young's inequality.
\end{proof}

\begin{lemma}[Dual space]
\label{lem:dual_space}
For $s \in \R$, $p,q \in [1, \infty)$ we have
\begin{equation*}
\left(M_{p,q}^s(\R^d) \right)^* = M_{p',q'}^{-s}(\R^d),
\end{equation*}
with the duality pairing given by
\begin{equation*}
\dup{v}{u}_{M_{p', q'}^{-s} \times M_{p, q}^s} =
\sum_{k \in \Z^d} \sum_{l \in M} \int_{\R^d}
\overline{\Box_{k + l} v} \Box_k u \D{x}
\end{equation*}
for any $v \in M_{p', q'}^{-s}(\R^d)$ and any $u \in M_{p, q}^s(\R^d)$,
and where $M$ is as in the proof of Lemma \ref{lem:bilinear}
(c.f. \cite[Theorem 3.1]{wang2007}).
\end{lemma}

\begin{theorem}[Complex interpolation]
\label{thm:complex_interpolation}
For $p_1, q_1 \in [1, \infty)$, $p_2, q_2 \in [1, \infty]$,
$s_1, s_2 \in \R$ and $\theta \in (0,1)$ one has
\begin{equation*}
\left[M_{p_1, q_1}^{s_1}(\R^d), M_{p_2, q_2}^{s_2}(\R^d)\right]_\theta =
M_{p,q}^s(\R^d),
\end{equation*}
with
\begin{equation*}
\frac{1}{p} = \frac{1 - \theta}{p_1} + \frac{\theta}{p_2}, \quad
\frac{1}{q} = \frac{1 - \theta}{q_1} + \frac{\theta}{q_2}, \quad
s = (1 - \theta) s_1 + \theta s_2
\end{equation*}
(see \cite[Theorem 6.1 (D)]{feichtinger1983}).
\end{theorem}

We are now ready to state and prove the following.
\begin{theorem}[Schr\"odinger propagator bound]
\label{thm:schroedinger}
There is a constant $C > 0$ such that for any $d \in \N$,
$p, q \in [1, \infty]$ and $s \in \R$ the inequality
\begin{equation}
\label{eqn:schroedinger_bound}
\norm{e^{\iu t \Delta}}_{\BLOBS(M_{p,q}^s(\R^d))} \leq
C^d (1 + \abs{t})^{d \abs{\frac{1}{2} - \frac{1}{p}}}
\end{equation}
holds for all $t \in \R$. Furthermore, the exponent of the time
dependence is sharp.
\end{theorem}
The boundedness has been obtained e.g. in
\cite[Theorem 1]{benyi2007} whereas the sharpness was proven in
\cite[Proposition 4.1]{cordero2009}. If $q < \infty$, then
$(e^{\iu t \Delta})_{t \in \R}$ is a $C_0$-group on $M_{p, q}^s$, i.e.,
\begin{equation*}
\lim_{t \to 0} \norm{e^{\iu t \Delta}f - f}_{M_{p, q}^s} = 0 \qquad
\forall f \in M_{p, q}^s
\end{equation*}
(see e.g. \cite[Proposition 3.5]{chaichenets2018}). This is not true for
$q = \infty$ and we refer to \cite{kunstmann2019} for this more subtle
case.

\begin{proof}[Proof of Theorem \ref{thm:schroedinger}]
By definition, we have
\begin{equation*}
(V_g e^{\iu t \Delta} f)(x, \xi) =
e^{-\iu t \abs{\xi}^2} (V_{e^{\iu t \Delta} g} f)(x + 2 t \xi, \xi)
\end{equation*}
for any $f \in \LSS'(\R^d)$, any $(x, \xi) \in \R^d \times \R^d$,
and any $t \in \R$, i.e. the Schr\"odinger time evolution of the initial
data can be interpreted as the time evolution of the window function.
The price for changing from window $g_0$ to window $g_1$ is
$\norm{V_{g_0} g_1}_{L^1(\R^d \times \R^d)}$ by
\cite[Proposition 11.3.2 (c)]{groechenig2001}. For
$g(x) = e^{-\abs{x}^2}$ one explicitly calculates
\begin{equation*}
\norm{V_{e^{- \iu t \Delta} g} g}_{L^1(\R^d \times \R^d)} = C^d
\left(1 + \abs{t} \right)^{\frac{d}{2}},
\end{equation*}
which proves the claimed bound for $p \in \set{1, \infty}$. Conservation
for $p = 2$ is easily seen from \eqref{eqn:equivalent_norm}.
Complex interpolation between the cases $p = 2$ and $p = \infty$ yields
\eqref{eqn:schroedinger_bound} for $p \in [2, \infty]$. The remaining
case $p \in (1, 2)$ is covered by duality.

Optimality in the case $p \in [1, 2]$ is proven by choosing
the window $g$ and the argument $f$ to be a Gaussian and explicitly
calculating $\norm{e^{\iu t \Delta} f}_{M_{p,q}^s} \approx
\left(1 + \abs{t}\right)^{d\left(\frac{1}{p} - \frac{1}{2}\right)}$.
This implies the optimality for $p \in (2, \infty]$ by duality.
\end{proof}

\section{Littlewood-Paley theory}
\label{sec:littlewood_paley}
In this section we extend some ideas of the Littlewood-Paley
decomposition from Sobolev spaces $H^s(\R^d)$ to modulation spaces
$M_{p,q}^s(\R^d)$. The inspiration for this was
\cite[Chapter II]{alinhac2007}.

Observe, that for any $\xi \in \R^d$ one has
\begin{eqnarray*}
\sum_{l = 0}^\infty \phi_l (\xi) =
\phi_0(\xi) + 
\lim_{N \to \infty} \sum_{l=1}^N
\left[\phi_1\left(\frac{\xi}{2^l}\right) -
\phi_1\left(\frac{\xi}{2^{l-1}}\right)\right] = 
\lim_{N \to \infty} \phi_0\left(\frac{\xi}{2^N}\right) = 1,
\end{eqnarray*}
i.e. $\set{\phi_0, \phi_1, \phi_2, \ldots}$ is a smooth partition of
unity. Moreover, $\supp(\phi_l) \subseteq A_l$ for any $l \in \N_0$,
where
\begin{equation*}
A_0 \coloneqq \set{\xi \in \R^d| \, \abs{\xi} \leq 1}
\quad \text{and} \quad
A_l \coloneqq \set{\xi \in \R^d| \, 2^{l - 2} \leq \abs{\xi} \leq 2^l}
\qquad \forall l \in \N.
\end{equation*}
The symbols of the dyadic decomposition operators satisfy
\begin{equation*}
\norm{\hat{\phi}_l}_1 =
\norm{\FT \left[\phi_1\left(\frac{\cdot}{2^{l-1}} \right)\right]}_1 =
\norm{2^{l-1} \hat{\phi}_1(2^{l-1} \cdot )}_1 =
\norm{\hat{\phi}_1}_1 \leq
2 \norm{\hat{\phi}_0}_1
\end{equation*}
for all $l \in \N$. Applying Lemma \ref{lem:bernstein} shows that
for any $l \in \N_0$ and any $f \in \LSS'(\R^d)$ one has that
$\Delta_l f \in C^\infty$ and any of its derivates has at most
polynomial growth. Furthermore, $\norm{\Delta_l}_{\BLOBS(L^p(\R^d))}$ is
bounded independently of $l \in \N_0$ and $p \in [1, \infty]$.

\begin{proof}[Proof of Theorem \ref{thm:littlewood_paley}]
We start by gathering some useful facts. Fix $l \in \N_0$ and
$k \in \Z^d$. Recall, that $\supp(\phi_l) \subseteq A_l$ and
$\supp(\sigma_k) \subseteq B_{\sqrt{d}}(k)$. Hence,
\begin{equation}
\label{eqn:supports}
\Box_k \Delta_l \not \equiv 0 \Rightarrow
k \in A'_l \coloneqq \set{k' \in \Z^d| \,
2^{l-2} - \sqrt{d} \leq \abs{k'} \leq 2^l + \sqrt{d}}.
\end{equation}

On $A'_l$ the Japanese bracket can be controlled. In fact, for all
$t \in \R$ we have
\begin{equation}
\label{eqn:jb_control}
\jb{k}^t \approx 2^{lt},
\end{equation}
where the implicit constant does not depend on $l$.

Finally, observe that $k \in A'_l$ is satisfied for only finitely many
$l \in \N_0$, whose number is independent of $k \in \Z^d$, i.e.
\begin{equation}
\label{eqn:intersections}
\sum_{l=0}^\infty \Ind_{A'_l}(k) \lesssim 1,
\end{equation}
where the implicit constant depends on $d$ only.

\begin{itemize}
\item $\gtrsim$: Consider $q < \infty$ first. By
\eqref{eqn:equivalent_norm}, \eqref{eqn:supports}, Bernstein multiplier
estimate, \eqref{eqn:jb_control} and \eqref{eqn:intersections} we have
\begin{eqnarray*}
& &
\norm{\left(2^{ls} \norm{\Delta_l f}_{M_{p,q}}\right)_l}_q \\
& \approx &
\left(\sum_{l=0}^\infty 2^{lsq}
\sum_{k \in \Z^d} \norm{\Box_k \Delta_l f}_p^q \right)^{\frac{1}{q}}
\lesssim
\left(\sum_{l=0}^\infty \sum_{k \in A'_l}
2^{lsq} \norm{\Box_k f}_p^q \right)^{\frac{1}{q}} \\
& \approx & 
\left(\sum_{l=0}^\infty \sum_{k \in \Z^d} \Ind_{A'_l} (k)
\jb{k}^{qs} \norm{\Box_k f}_p^q \right)^{\frac{1}{q}} \lesssim
\norm{f}_{M_{p,q}^s}.
\end{eqnarray*}
Similarly, for $q = \infty$, we have
\begin{eqnarray*}
\norm{\left(2^{ls} \norm{\Delta_l f}_{M_{p,\infty}}\right)_l}_\infty
& = &
\sup_{l \in \N_0} 2^{ls} \sup_{k \in \Z^d} \norm{\Box_k \Delta_l f}_p
\\
& \lesssim &
\sup_{l \in \N_0} \sup_{k \in A'_l} \jb{k}^s \norm{\Box_k f}_p \approx
\norm{f}_{M_{p,\infty}^s}.
\end{eqnarray*}

\item $\lesssim$: Again, consider $q < \infty$ first. By
\eqref{eqn:equivalent_norm}, $f = \sum_{l = 0}^\infty \Delta_l f$ in
$\LSS'$ and \eqref{eqn:supports} we have
\begin{eqnarray*}
\norm{f}_{M_{p,q}^s} & \lesssim &
\left(\sum_{k \in \Z^d} \jb{k}^{qs} \left(\sum_{l=0}^\infty
\norm{\Box_k \Delta_l f}_p\right)^q \right)^{\frac{1}{q}} \\
& \lesssim &
\left(\sum_{k \in \Z^d} \jb{k}^{qs} \left(\sum_{l=0}^\infty
\Ind_{A'_l}(k) \norm{\Box_k \Delta_l f}_p\right)^q \right)^{\frac{1}{q}}.
\end{eqnarray*}
For each $k \in \Z^d$ the sum over $l$ contains only finitely many
non-vanishing summands and their number is independent of $k$ by
\eqref{eqn:intersections}. H\"older's inequality estimates the last term
against
\begin{eqnarray*}
\left(\sum_{k \in \Z^d} \jb{k}^{qs} \sum_{l=0}^\infty
\Ind_{A'_l}(k) \norm{\Box_k \Delta_l f}_p^q \right)^{\frac{1}{q}}
& \approx &
\left(\sum_{l=0}^\infty 2^{lsq} \sum_{k \in \Z^d}
\Ind_{A'_l}(k) 
\norm{\Delta_l \Box_k f}_p^q \right)^{\frac{1}{q}} \\
& \leq &
\norm{\left(2^{ls} \norm{\Delta_l f}_{M_{p,q}} \right)_l}_q,
\end{eqnarray*}
where we additionally used \eqref{eqn:jb_control}. The proof for
$q = \infty$ is along the same lines.
\end{itemize}
\end{proof}

The individual parts of the Littlewood-Paley decomposition had their
Fourier transform supported in almost disjoint dyadic annuli. Theorem
\ref{thm:littlewood_paley} characterized elements of modulation spaces
by the decay of these parts. The following lemma provides a sufficient
condition for the case of non-disjoint balls.
\begin{lemma}[Sufficient condition]
\label{lem:sufficient}
Let $1 \leq q \leq \infty$ and $s > 0$. For $m \in \N_0$ let
$f_m \in \LSS'(\R^d)$ be such that
\begin{equation*}
\supp(\hat{f}_m) \subseteq B_m \coloneqq
\set{\xi \in \R^d \big| \abs{\xi} \leq 2^m} \qquad \forall m \in \N_0.
\end{equation*}
Set $f \coloneqq \sum_{m=0}^\infty f_m$ in $\LSS'(\R^d)$. Then
\begin{equation*}
\norm{f}_{M_{p,q}^s(\R^d)} \lesssim
\norm{\left(2^{ms} \norm{f_m}_{M_{p,q}(\R^d)} \right)_{m \in \N_0}}_q,
\end{equation*}
where the implicit constant depends on $d$ and $s$ only.
\end{lemma}

\begin{proof}
Observe, that $A_l \cap B_m = \emptyset$ if $l > m + 2$. Hence, we have
\begin{eqnarray*}
\norm{f}_{M_{p,q}^s} & \approx &
\norm{\left(2^{ls} \norm{\Delta_l f}_{M_{p,q}}\right)_l}_q
\lesssim
\norm{\left(2^{ls} \sum_{m = l}^\infty
\norm{\Delta_l f_m}_{M_{p,q}}\right)_l}_q \\
& \lesssim &
\norm{\left(2^{ls} \sum_{m = l}^\infty
\norm{f_m}_{M_{p,q}}\right)_l}_q,
\end{eqnarray*}
where we additionally used Theorem \ref{thm:littlewood_paley} and
Bernstein multiplier estimate. From now on, we assume
$q \in (1, \infty)$, as the proof for the other cases is easier and
follows the same lines. H\"older's inequality and geometric sum formula
estimates the last term against
\begin{eqnarray*}
& &
\left(\sum_{l=0}^\infty \left(\sum_{m = l}^\infty
2^{ls} \norm{f_m}_{M_{p,q}}\right)^q \right)^{\frac{1}{q}} \\
& = &
\left(\sum_{l=0}^\infty \left(\sum_{m = l}^\infty
2^{\frac{(l-m)s}{q'}} \times 2^{\frac{(l-m)s}{q}} 2^{ms}
\norm{f_m}_{M_{p,q}}\right)^q
\right)^{\frac{1}{q}} \\
& \leq &
\left(\sum_{l=0}^\infty
\left(\sum_{m = l}^\infty 2^{(l-m)s}\right)^{\frac{q}{q'}}
\left(\sum_{m = l}^\infty 2^{(l-m)s} 2^{msq}
\norm{f_m}_{M_{p,q}}^q\right)\right)^{\frac{1}{q}} \\
& \approx &
\left(\sum_{m=0}^\infty
\sum_{l = 0}^m 2^{(l-m)s} 2^{msq}
\norm{f_m}_{M_{p,q}}^q\right)^{\frac{1}{q}} \\
& \approx &
\norm{\left(2^{ms} \norm{f_m}_{M_{p,q}}\right)_{m}}_q,
\end{eqnarray*}
finishing the proof.
\end{proof}

\section{Algebra property and H\"older-type inequality}
\label{sec:algebra}
Main goal of this section is to prove Theorem
\ref{thm:algebra}, which was inspired by the fact that
$H^s(\R^d) \cap L^\infty(\R^d)$ is a Banach *-algebra with respect to
pointwise multiplication for $s \geq 0$.
\begin{proof}[Proof of Theorem \ref{thm:algebra}]
Parts \ref{it:regularity_fourier_integrability} and
\ref{it:contembedding} of Lemma \ref{lem:embeddings} prove the claimed
embedding. Continuity of complex conjugation is obvious from
\eqref{eqn:equivalent_norm}. Continuity of multiplication follows by the
paraproduct argument
\begin{equation*}
f g =
\left(\sum_{l = 0}^\infty \Delta_l f\right)
\left(\sum_{m = 0}^\infty \Delta_m g\right) =
\sum_{l = 0}^\infty \underbrace{
\left( \Delta_l f
\sum_{m = 0}^l \Delta_m g\right)}_{=: u_l} +
\sum_{m = 1}^\infty \underbrace{
\left(\Delta_m g
\sum_{l = 0}^{m-1} \Delta_l f\right)}_{=: v_m}.
\end{equation*}
Observe, that for any $l, m \in \N_0$ we have
$\supp (\hat{u}_l) \subseteq B_{l+1}$ and
$\supp (\hat{v}_m) \subseteq B_m$ by the properties of convolution.
Hence, Lemma \ref{lem:sufficient} could be applied.
Bilinear estimate from Lemma \ref{lem:bilinear} and Theorem
\ref{thm:littlewood_paley} show
\begin{eqnarray*}
\norm{\left(2^{ls} \norm{u_l}_{M_{p,q}}\right)_l}_q & \leq &
\norm{\left(2^{ls} \norm{\Delta_l f}_{M_{p,q}}\right)_l}_q
\sum_{m=0}^\infty \norm{\Delta_m g}_{M_{\infty,1}} \approx
\norm{f}_{M_{p,q}^s} \norm{g}_{M_{\infty, 1}}.
\end{eqnarray*}
The same argument yields $\norm{\sum_{m=1}^\infty v_m}_{M_{p,q}^s}
\lesssim \norm{f}_{M_{\infty,1}}\norm{g}_{M_{p,q}^s}$ and finishes the
proof.
\end{proof}

The analogon of Theorem \ref{thm:algebra} for sequence spaces
is stated in
\begin{lemma}[Algebra property]
\label{lem:algebra_property_seq}
Let $1 \leq q \leq \infty$ and $s \geq 0$. Then
$l_s^q(\Z^d) \cap l^1(\Z^d)$ is a Banach algebra with respect to
convolution
\begin{equation}
\label{eqn:discrete_convolution}
(a_{l}) \ast (b_{m}) =
\left(\sum_{m \in \Z^d }
a_{k - m} b_{m} \right)_{k \in \Z^d},
\end{equation}
which is well-defined, as the series above always converge absolutely.

Furthermore, if $q > 1$ and $s > d \left(1 - \frac{1}{q}\right)$ or
$q = 1$, then $l^q_s(\Z^d) \hookrightarrow l^1(\Z^d)$, so in particular
$l^q_s(\Z^d)$ is a Banach algebra, in that case.
\end{lemma}
This result is not new, see e.g. \cite[Satz 3.7]{feichtinger1979}. The
proof given there works in much more general situations and
relies on the fact that for $s \geq 0$ the weight function
$x \mapsto \jb{x}^s$ is \emph{weakly sub-additive} (a notion going back
to \cite{brandenburg1975}, at least), i.e.,
\begin{equation*}
\jb{x + y}^s \lesssim \jb{x}^s + \jb{y}^s \qquad \forall x, y \in \R.
\end{equation*} 
For the sake of a self-contained presentation, let us remark that
another proof can be given using the same techniques as for the proof
of Theorem \ref{thm:algebra}, i.e. by proving analoga of Theorem
\ref{thm:littlewood_paley} and Lemma \ref{lem:sufficient} for the
weighted sequence spaces. Yet another approach is to notice that by
definition
\begin{equation*}
\norm{\sum_{k \in \Z^d} a_k e^{\iu k x}}_{M_{\infty, q}^s} \approx
\norm{(a_k)}_{l^q_s}
\end{equation*}
and hence, by Theorem \ref{thm:algebra}, one has
\begin{eqnarray*}
& &
\norm{(a_k) \ast (b_k)}_{l^q_s} \\
& \approx &
\norm{\left(\sum_{k \in \Z^d} a_k e^{\iu k x}\right) \cdot
\left(\sum_{k \in \Z^d} b_k e^{\iu k x}\right)}_{M_{\infty, q}^s}
\\
& \lesssim &
\norm{\sum_{k \in \Z^d} a_k e^{\iu k x}}_{M_{\infty, q}^s}
\norm{\sum_{k \in \Z^d} b_k e^{\iu k x}}_{M_{\infty, 1}} +
\norm{\sum_{k \in \Z^d} a_k e^{\iu k x}}_{M_{\infty, 1}}
\norm{\sum_{k \in \Z^d} b_k e^{\iu k x}}_{M_{\infty, q}^s} \\
& \approx &
\norm{(a_k)}_{l^q_s} \norm{(b_k)}_{l^1} + 
\norm{(a_k)}_{l^1} \norm{(b_k)}_{l^q_s}.
\end{eqnarray*}

We are now ready to give a
\begin{proof}[Proof of Theorem \ref{thm:hoelder}]
We arrive, as for equation \eqref{eqn:gen_hoelder} in the proof of Lemma
\ref{lem:bilinear}, at
\begin{equation*}
\norm{\Box_{k}(fg)}_p \lesssim
\sum_{m \in M} \left(
\left(\norm{\Box_{l}(f)}_{p_1}\right)_{l} \ast
\left(\norm{\Box_{n + m}(g)}_{p_2}\right)_n \right)(k)
\end{equation*}
pointwise in $k \in \Z^d$. By the algebra property from Lemma
\ref{lem:algebra_property_seq}, it follows that
\begin{equation*}
\norm{fg}_{M_{p,q}^s} \lesssim
\norm{\left(\norm{\Box_{l} f}_{p_1} \right)_{l}}_{q,s}
\left(\sum_{m \in M}
\norm{
\left(\norm{\Box_{n + m} g}_{p_2} \right)_{n}}_{q,s}
\right)
\end{equation*}
and the first factor is already $\norm{f}_{M_{p,q}^s}$. Finally, we
remove the sum over $m$ in the second factor
\begin{equation*}
\sum_{m \in M}
\norm{
\left(\norm{\Box_{n + m} g}_{p_2} \right)_{n}}_{q,s}
\lesssim
\norm{g}_{M_{p_2,q}^s}
\end{equation*}
applying Peetre's inequality $\jb{k + l}^s \leq
2^{\abs{s}} \jb{k}^s \jb{l}^{\abs{s}}$ (see e.g.
\cite[Proposition 3.3.31]{ruzhansky2010}).

Let us finish the proof remarking that the only estimate involving
``$p$''s we used was H\"older's inequality and thus the implicit constant
indeed does not depend on $p$, $p_1$ or $p_2$.
\end{proof}

\section{Proof of the local well-posedness, Theorem \ref{thm:lwp}.}
\label{sec:local_well-posedness}
Theorem \ref{thm:algebra} immediately implies that $X(T)$ is
a Banach *-algebra, i.e.
\begin{eqnarray*}
\norm{uv}_{X(T)} & = &
\sup_{0 \leq t \leq T} \norm{uv(\cdot, t)}_X \lesssim
\left(\sup_{0 \leq s \leq T} \norm{u(\cdot, s)}_X \right)
\left(\sup_{0 \leq t \leq T} \norm{v(\cdot, t)}_X \right) \\
& = &
\norm{u}_{X(T)} \norm{v}_{X(T)}.
\end{eqnarray*}
For $R > 0$ we denote by $M(R,T) =
\set{u \in X(T) \Big| \norm{u}_{X(T)} \leq R}$ the closed ball of radius
$R$ in $X(T)$ centered at the origin.
We show that for some $R, T > 0$ the right-hand side of
\eqref{eqn:duhamel},
\begin{equation}
\label{eqn:contraction}
\left(\mathcal{T} u\right)(\cdot,t) \coloneqq
e^{\iu t \Delta} u_0 \pm \iu
\int_{0}^t e^{\iu(t-\tau) \Delta}
\left(\abs{u}^2 u(\cdot, \tau)\right) \D{\tau} \qquad
(\forall t \in [0, T]),
\end{equation}
defines a contractive self-mapping $\mathcal{T} = \mathcal{T}(u_0):
M_{R,T} \to M_{R,T}$.

To that end, let us observe that Theorem \ref{thm:schroedinger}
implies the \emph{homogeneous estimate}
\begin{equation*}
\norm{t \mapsto e^{\iu t \Delta} v}_X \leq
C_0 (1 + T)^{\frac{d}{2}} \norm{v}_X \qquad
(\forall v \in X),
\end{equation*}
which, together with the algebra property of $X(T)$, proves the
\emph{inhomogeneous estimate}
\begin{eqnarray*}
& &
\norm{
\int_0^t e^{\iu (t - \tau) \Delta}
\left(\abs{u}^2 u(\cdot, \tau) \right) \D{\tau}}_X \\
& \leq &
C_0 (1 + T)^{\frac{d}{2}} \int_0^t
\norm{\abs{u}^2 u(\cdot, \tau)}_X \D{\tau} \leq
C_0 C_1 T (1 + T)^{\frac{d}{2}} \norm{u}_X^3,
\end{eqnarray*}
holding for $0 \leq t \leq T$ and $u \in X(T)$. 

Applying the triangle inequality in \eqref{eqn:contraction} yields
\begin{equation*}
\norm{\mathcal{T} u}_X \leq
C_0 (1 + T)^{\frac{d}{2}} (\norm{u_0}_X + C_1 T R^3)
\end{equation*}
for any $u \in M(R,T).$ Thus, $\mathcal{T}$ maps $M(R,T)$ into itself
for $R = 2 C_0 C_1 \norm{u_0}_X$ and $T$ small enough. Furthermore,
\begin{equation*}
\abs{u}^2 u - \abs{v}^2 v =
(u - v) \abs{u}^2 + (\overline{u}u - \overline{v}v)v =
(u - v)(\abs{u}^2 + \overline{u} v) +
(\overline{u} - \overline{v})v^2
\end{equation*}
and hence
\begin{equation*}
\norm{\mathcal{T}u - \mathcal{T}v}_{X(T)} \lesssim
T (1 + T)^{\frac{d}{2}} R^2 \norm{u - v}_{X(T)}
\end{equation*}
for $u, v \in M(R,T)$, where we additionally used the algebra property
of $X(T)$ and the homogeneous estimate. Taking $T$ sufficiently small
makes $\mathcal{T}$ a contraction.

Banach's fixed-point theorem implies the existence and uniqueness of a
mild solution up to the \emph{guaranteed time of existence} $T_0 =
T_0\left(\norm{u_0}_X \right) \approx \norm{u_0}_X^{-2} > 0$. Uniqueness
of the maximal solution and the blow-up alternative now follow easily by
the usual contradiction argument.

For the proof of the Lipschitz continuity, let us notice that for any
$r > \norm{u_0}_X$, $v_0 \in B_r(0)$ and $0 < T \leq T_0(r)$
we have
\begin{eqnarray*}
\norm{u - v}_{X(T)} & = &
\norm{\mathcal{T}(u_0) u - \mathcal{T}(v_0)v}_{X(T)} \\
& \lesssim &
(1 + T)^{\frac{d}{2}} \norm{u_0 - v_0}_X +
T (1 + T)^{\frac{d}{2}} R^2 \norm{u - v}_{X(T)},
\end{eqnarray*}
where $v$ is the mild solution corresponding to the initial data $v_0$
and $R = 2Cr$, similar to the above. Collecting terms containing
$\norm{u - v}_{X(T)}$ shows Lipschitz continuity with constant
$L = L(r)$ for sufficiently small $T$, say $T_l = T_l(r)$. For arbitrary
$0 < T' < T_*$ put $r = 2 \norm{u}_{X(T')}$ and divide $[0,T']$ into $n$
subintervals of length $\leq T_l$. The claim follows for
$V = B_{\delta} (u_0)$ where $\delta =
\frac{\norm{u_0}_X}{L^n}$ by iteration. This concludes the proof.
\hfill\ensuremath{\square}

\section*{Acknowledgments}
The authors would like to thank Professor H.G. Feichtinger for comments
and remarks on an earlier version of this paper.

We gratefully acknowledge financial support by the Deutsche
Forschungsgemeinschaft (DFG) through CRC 1173.

\bibliography{literature.bib}

\begin{thebibliography}{WHHG11}
  \providebibliographyfont{name}{}%
  \providebibliographyfont{lastname}{}%
  \providebibliographyfont{title}{\emph}%
  \providebibliographyfont{jtitle}{\btxtitlefont}%
  \providebibliographyfont{etal}{\emph}%
  \providebibliographyfont{journal}{}%
  \providebibliographyfont{volume}{}%
  \providebibliographyfont{ISBN}{\MakeUppercase}%
  \providebibliographyfont{ISSN}{\MakeUppercase}%
  \providebibliographyfont{url}{\url}%
  \providebibliographyfont{numeral}{}%
  \expandafter\btxselectlanguage\expandafter {\btxfallbacklanguage}

\btxselectlanguage {english}
\bibitem [{AG07}]{alinhac2007}
\btxnamefont {Serge \btxlastnamefont {Alinhac}} \btxandlong {}\ \btxnamefont
  {Patrick \btxlastnamefont {G\'erard}}\btxauthorcolon\ \btxtitlefont
  {Pseudo-differential Operators and the {Nash-Moser} Theorem}, \btxvolumelong
  {}~\btxvolumefont {82} \btxofserieslong {}\ \btxtitlefont {Graduate Studies
  in Mathematics}.
\newblock \btxpublisherfont {American Mathematical Society}, Providence, Rhode
  Island, 2007\ifbtxprintISBN {, \mbox{\btxISBN~\btxISBNfont
  {978-0-8218-3454-1}}}.

\bibitem [{BGOR07}]{benyi2007}
\btxnamefont {\'Arp\'ad \btxlastnamefont {B\'enyi}}, \btxnamefont {Karlheinz
  \btxlastnamefont {Gr\"ochenig}}, \btxnamefont {Kasso\btxfnamespacelong
  Akochay\'e \btxlastnamefont {Okoudjou}}\btxandcomma {} \btxandlong {}\
  \btxnamefont {Luke\btxfnamespacelong Gervase \btxlastnamefont
  {Rogers}}\btxauthorcolon\ \btxjtitlefont {\btxifchangecase {Unimodular
  {Fourier} multipliers for modulation spaces}{Unimodular {Fourier} multipliers
  for modulation spaces}}.
\newblock \btxjournalfont {Journal of Functional Analysis}, 246(2):366--384,
  2007\ifbtxprintISSN {, \mbox{\btxISSN~\btxISSNfont {0022-1236}}}.
\newblock {\latintext \btxurlfont{https://doi.org/10.1016/j.jfa.2006.12.019}}.

\bibitem [{BO09}]{benyi2009}
\btxnamefont {\'Arp\'ad \btxlastnamefont {B\'enyi}} \btxandlong {}\
  \btxnamefont {Kasso\btxfnamespacelong Akochay\'e \btxlastnamefont
  {Okoudjou}}\btxauthorcolon\ \btxjtitlefont {\btxifchangecase {Local
  well-posedness of nonlinear dispersive equations on modulation spaces}{Local
  well-posedness of nonlinear dispersive equations on modulation spaces}}.
\newblock \btxjournalfont {Bulletin of the London Mathematical Society},
  41(3):549--558, 2009\ifbtxprintISSN {, \mbox{\btxISSN~\btxISSNfont
  {0024-6093}}}.
\newblock {\latintext \btxurlfont{https://doi.org/10.1112/blms/bdp027}}.

\bibitem [{Bra75}]{brandenburg1975}
\btxnamefont {L.\btxfnamespacelong H. \btxlastnamefont
  {Brandenburg}}\btxauthorcolon\ \btxjtitlefont {\btxifchangecase {On
  identifying the maximal ideals in {B}anach algebras}{On Identifying the
  Maximal Ideals in {B}anach Algebras}}.
\newblock \btxjournalfont {Journal of Mathematical Analysis and Applications},
  50:489 -- 510, 1975\ifbtxprintISSN {, \mbox{\btxISSN~\btxISSNfont
  {0022-247X}}}.
\newblock {\latintext
  \btxurlfont{https://doi.org/10.1016/0022-247X(75)90006-2}}.

\bibitem [{Cha18}]{chaichenets2018}
\btxnamefont {Leonid \btxlastnamefont {Chaichenets}}\btxauthorcolon\
  \btxtitlefont {Modulation spaces and nonlinear {S}chr\"odinger equations}.
\newblock \btxphdthesis {}, Karlsruhe Institute of Technology (KIT), 2018.
\newblock {\latintext \btxurlfont{https://doi.org/10.5445/IR/1000088173}}.

\bibitem [{CHKP19}]{pattakos2018a}
\btxnamefont {Leonid \btxlastnamefont {Chaichenets}}, \btxnamefont {Dirk
  \btxlastnamefont {Hundertmark}}, \btxnamefont {Peer \btxlastnamefont
  {Kunstmann}}\btxandcomma {} \btxandlong {}\ \btxnamefont {Nikolaos
  \btxlastnamefont {Pattakos}}\btxauthorcolon\ \btxjtitlefont {\btxifchangecase
  {Nonlinear {S}chr\"odinger equation, differentiation by parts and modulation
  spaces}{Nonlinear {S}chr\"odinger equation, differentiation by parts and
  modulation spaces}}.
\newblock \btxjournalfont {Journal of Evolution Equations}, 2019\ifbtxprintISSN
  {, \mbox{\btxISSN~\btxISSNfont {1424-3202}}}.
\newblock {\latintext \btxurlfont{https://doi.org/10.1007/s00028-019-00501-z}}.

\bibitem [{CN09}]{cordero2009}
\btxnamefont {Elena \btxlastnamefont {Cordero}} \btxandlong {}\ \btxnamefont
  {Fabio \btxlastnamefont {Nicola}}\btxauthorcolon\ \btxjtitlefont
  {\btxifchangecase {Sharpness of some properties of {Wiener} amalgam and
  modulation spaces}{Sharpness of some properties of {Wiener} amalgam and
  modulation spaces}}.
\newblock \btxjournalfont {Bulletin of the Australian Mathematical Society},
  80(1):105 -- 116, 2009\ifbtxprintISSN {, \mbox{\btxISSN~\btxISSNfont
  {0004-9727}}}.
\newblock {\latintext \btxurlfont{https://doi.org/10.1017/S0004972709000070}}.

\btxselectlanguage {ngerman}
\bibitem [{Fei79}]{feichtinger1979}
\btxnamefont {Hans\btxfnamespacelong Georg \btxlastnamefont
  {Feichtinger}}\btxauthorcolon\ \btxjtitlefont {\btxifchangecase
  {Gewichtsfunktionen auf lokalkompakten {G}ruppen}{Gewichtsfunktionen auf
  lokalkompakten {G}ruppen}}.
\newblock \btxjournalfont {Sitzungsberichte der {\"O}sterreicherischen
  {A}kademie der {W}issenschaften}, 188(8 -- 10):451 -- 471,
  1979\ifbtxprintISSN {, \mbox{\btxISSN~\btxISSNfont {0029-8816}}}.

\btxselectlanguage {english}
\bibitem [{Fei80}]{feichtinger1980}
\btxnamefont {Hans\btxfnamespacelong Georg \btxlastnamefont
  {Feichtinger}}\btxauthorcolon\ \btxjtitlefont {\btxifchangecase {Banach
  convolution algebras of {Wiener} type}{Banach convolution algebras of
  {Wiener} type}}.
\newblock \btxjournalfont {Functions, Series, Operators}, 35:509 -- 524, 1980.
\newblock {\latintext
  \btxurlfont{https://www.univie.ac.at/nuhag-php/bibtex/open_files/fe83_wientyp1.pdf}}.

\bibitem [{Fei83}]{feichtinger1983}
\btxnamefont {Hans\btxfnamespacelong Georg \btxlastnamefont
  {Feichtinger}}\btxauthorcolon\ \btxtitlefont {Modulation spaces on locally
  compact abelian groups}.
\newblock \btxpublisherfont {University Vienna}, 1983.
\newblock {\latintext
  \btxurlfont{https://www.univie.ac.at/nuhag-php/bibtex/open_files/120_ModICWA.pdf}}.

\expandafter\btxselectlanguage\expandafter {\btxfallbacklanguage}
\bibitem [{GFWZ18}]{guo2018}
\btxnamefont {Weichao \btxlastnamefont {Guo}}, \btxnamefont {Dashan
  \btxlastnamefont {Fan}}, \btxnamefont {Huoxiong \btxlastnamefont
  {Wu}}\btxandcomma {} \btxandlong {}\ \btxnamefont {Guoping \btxlastnamefont
  {Zhao}}\btxauthorcolon\ \btxjtitlefont {\btxifchangecase {Sharp weighted
  convolution inequalities and some applications}{Sharp weighted convolution
  inequalities and some applications}}.
\newblock \btxjournalfont {Studia Mathematica}, 241(3):201 -- 239,
  2018\ifbtxprintISSN {, \mbox{\btxISSN~\btxISSNfont {0039-3223}}}.
\newblock {\latintext \btxurlfont{https://doi.org/10.4064/sm8583-5-2017}}.

\btxselectlanguage {english}
\bibitem [{Gr{\"o}01}]{groechenig2001}
\btxnamefont {Karlheinz \btxlastnamefont {Gr{\"o}chenig}}\btxauthorcolon\
  \btxtitlefont {Foundations of time-frequency analysis}.
\newblock Applied and Numerical Harmonic Analysis. \btxpublisherfont
  {Birkh\"auser}, Boston, 2001\ifbtxprintISBN {, \mbox{\btxISBN~\btxISBNfont
  {978-0-8176-4022-4}}}.
\newblock {\latintext \btxurlfont{https://doi.org/10.1007/978-1-4612-0003-1}}.

\bibitem [{Guo16}]{guo2016}
\btxnamefont {Shaoming \btxlastnamefont {Guo}}\btxauthorcolon\ \btxjtitlefont
  {\btxifchangecase {On the {1D} cubic {NLS} in an almost critical space}{On
  the {1D} cubic {NLS} in an almost critical space}}.
\newblock \btxjournalfont {Journal of Fourier Analysis and Applications},
  23(1):91 -- 124, 2016\ifbtxprintISSN {, \mbox{\btxISSN~\btxISSNfont
  {1531-5851}}}.
\newblock {\latintext \btxurlfont{https://doi.org/10.1007/s00041-016-9464-z}}.

\bibitem [{Kun19}]{kunstmann2019}
\btxnamefont {Peer\btxfnamespacelong Christian \btxlastnamefont
  {Kunstmann}}\btxauthorcolon\ \btxjtitlefont {\btxifchangecase {Modulation
  type spaces for generators of polynomially bounded groups and {S}chr\"odinger
  equations}{Modulation type spaces for generators of polynomially bounded
  groups and {S}chr\"odinger equations}}.
\newblock \btxjournalfont {Semigroup Forum}, 2019\ifbtxprintISSN {,
  \mbox{\btxISSN~\btxISSNfont {1432-2137}}}.
\newblock {\latintext \btxurlfont{https://doi.org/10.1007/s00233-019-10016-1}}.

\bibitem [{Pat18}]{pattakos2018}
\btxnamefont {Nikolaos \btxlastnamefont {Pattakos}}\btxauthorcolon\
  \btxjtitlefont {\btxifchangecase {{NLS} in the modulation space {$M_{2,
  q}(\R)$}}{{NLS} in the Modulation Space {$M_{2, q}(\R)$}}}.
\newblock \btxjournalfont {Journal of Fourier Analysis and Applications},
  2018\ifbtxprintISSN {, \mbox{\btxISSN~\btxISSNfont {1531-5851}}}.
\newblock {\latintext \btxurlfont{https://doi.org/10.1007/s00041-018-09655-9}}.

\bibitem [{RSW12}]{ruzhansky2012}
\btxnamefont {Michael \btxlastnamefont {Ruzhansky}}, \btxnamefont {Mitsuru
  \btxlastnamefont {Sugimoto}}\btxandcomma {} \btxandlong {}\ \btxnamefont
  {Baoxiang \btxlastnamefont {Wang}}\btxauthorcolon\ \btxtitlefont
  {\btxifchangecase {Modulation spaces and nonlinear evolution
  equations}{Modulation spaces and nonlinear evolution equations}}.
\newblock \Btxinlong {}\ \btxtitlefont {Evolution equations of hyperbolic and
  Schr\"odinger type}, \btxvolumelong {}\ \btxvolumefont {301}, \btxpageslong
  {}\ 267--283. \btxpublisherfont {Springer}, Basel, 2012\ifbtxprintISBN {,
  \mbox{\btxISBN~\btxISBNfont {978-3-0348-0453-0}}}.
\newblock {\latintext
  \btxurlfont{https://doi.org/10.1007/978-3-0348-0454-7_14}}.

\bibitem [{RT10}]{ruzhansky2010}
\btxnamefont {Michael\btxfnamespacelong Vladimirovich \btxlastnamefont
  {Ruzhansky}} \btxandlong {}\ \btxnamefont {Ville \btxlastnamefont
  {Turunen}}\btxauthorcolon\ \btxtitlefont {Pseudo-Differential Operators and
  Symmetries}.
\newblock \Btxnumberlong {}~2 \btxinserieslong {}\ \btxtitlefont
  {Pseudo-Differential Operators}. \btxpublisherfont {Birkh\"auser}, Basel,
  2010\ifbtxprintISBN {, \mbox{\btxISBN~\btxISBNfont {978-3-7643-8513-2}}}.
\newblock {\latintext \btxurlfont{https://doi.org/10.1007/978-3-7643-8514-9}}.

\bibitem [{STW11}]{sugimoto2011}
\btxnamefont {Mitsuru \btxlastnamefont {Sugimoto}}, \btxnamefont {Naohito
  \btxlastnamefont {Tomita}}\btxandcomma {} \btxandlong {}\ \btxnamefont
  {Baoxiang \btxlastnamefont {Wang}}\btxauthorcolon\ \btxjtitlefont
  {\btxifchangecase {Remarks on nonlinear operations on modulation
  spaces}{Remarks on nonlinear operations on modulation spaces}}.
\newblock \btxjournalfont {Integral Transforms and Special Functions}, 22(4 --
  5):351 -- 358, 2011\ifbtxprintISSN {, \mbox{\btxISSN~\btxISSNfont
  {1065-2469}}}.
\newblock {\latintext
  \btxurlfont{https://doi.org/10.1080/10652469.2010.541054}}.

\bibitem [{TJPT15}]{toft2015}
\btxnamefont {Joachim \btxlastnamefont {Toft}}, \btxnamefont {Karoline
  \btxlastnamefont {Johansson}}, \btxnamefont {Stevan \btxlastnamefont
  {Pilipović}}\btxandcomma {} \btxandlong {}\ \btxnamefont {Nenad
  \btxlastnamefont {Teofanov}}\btxauthorcolon\ \btxjtitlefont {\btxifchangecase
  {Sharp convolution and multiplication estimates in weighted spaces}{Sharp
  convolution and multiplication estimates in weighted spaces}}.
\newblock \btxjournalfont {Analysis and Applications}, 13(5):457 -- 480,
  2015\ifbtxprintISSN {, \mbox{\btxISSN~\btxISSNfont {0219-5305}}}.
\newblock {\latintext \btxurlfont{https://doi.org/10.1142/S0219530514500523}}.

\bibitem [{WH07}]{wang2007}
\btxnamefont {Baoxiang \btxlastnamefont {Wang}} \btxandlong {}\ \btxnamefont
  {Henryk \btxlastnamefont {Hudzik}}\btxauthorcolon\ \btxjtitlefont
  {\btxifchangecase {The global {Cauchy} problem for the {NLS} and {NLKG} with
  small rough data}{The global {Cauchy} problem for the {NLS} and {NLKG} with
  small rough data}}.
\newblock \btxjournalfont {Journal of Differential Equations}, 232(1):36--73,
  2007\ifbtxprintISSN {, \mbox{\btxISSN~\btxISSNfont {0022-0396}}}.
\newblock {\latintext \btxurlfont{https://doi.org/10.1016/j.jde.2006.09.004}}.

\bibitem [{WHHG11}]{wang2011}
\btxnamefont {Baoxiang \btxlastnamefont {Wang}}, \btxnamefont {Zhaohui
  \btxlastnamefont {Huo}}, \btxnamefont {Chengchun \btxlastnamefont
  {Hao}}\btxandcomma {} \btxandlong {}\ \btxnamefont {Zinhua \btxlastnamefont
  {Guo}}\btxauthorcolon\ \btxtitlefont {Harmonic Analysis Method for Nonlinear
  Evolution Equations, {I}}.
\newblock \btxpublisherfont {World Scientific}, 2011\ifbtxprintISBN {,
  \mbox{\btxISBN~\btxISBNfont {978-981-4360-73-9}}}.
\newblock {\latintext \btxurlfont{https://doi.org/10.1142/8209}}.

\bibitem [{WZG06}]{wang2006}
\btxnamefont {Baoxiang \btxlastnamefont {Wang}}, \btxnamefont {Lifeng
  \btxlastnamefont {Zhao}}\btxandcomma {} \btxandlong {}\ \btxnamefont {Boling
  \btxlastnamefont {Guo}}\btxauthorcolon\ \btxjtitlefont {\btxifchangecase
  {Isometric decomposition operators, function spaces {$E^\lambda_{p, q}$} and
  applications to nonlinear evolution equations}{Isometric decomposition
  operators, function spaces {$E^\lambda_{p, q}$} and applications to nonlinear
  evolution equations}}.
\newblock \btxjournalfont {Journal of Functional Analysis}, 233(1):1--39,
  2006\ifbtxprintISSN {, \mbox{\btxISSN~\btxISSNfont {0022-1236}}}.
\newblock {\latintext \btxurlfont{https://doi.org/10.1016/j.jfa.2005.06.018}}.

\end{thebibliography}
\bibliographystyle{babalpha-fl}

\end{document}